\documentclass[reqno,11pt]{article}
\usepackage{amsmath,amssymb}
\usepackage{amsthm,amscd,amsfonts,}
\usepackage{amsmath}
\usepackage{amsthm}
\usepackage{graphicx}
\usepackage{amsmath,amsthm}
\usepackage{amsmath,amstext,amsfonts}
\usepackage{subfigure}

\newtheorem{theorem}{Theorem}[section]
\newtheorem{lemma}{Lemma}[section]
\newtheorem{remark}{Remark}[section]

\newcommand{\R}{\mathbb{R}}
\numberwithin{equation}{section}

\begin{document}

\title{\textbf{A Global Attractivity in a Nonmonotone  Age-Structured Model with Age Dependent Diffusion and Death Rates}}
\author{M. Al-Jararha\thanks{%
Department of Mathematics, Yarmouk University, Irbid, Jordan, 21163. \texttt{%
mohammad.ja@yu.edu.jo}.}}
\date{}
\maketitle

\begin{abstract}
In this paper, we investigated the global attractivity of the  positive  constant steady state solution of the mature population $w(t,x)$  governed by  the age-structured model: 
\begin{equation*}
\left\{\begin{array}{ll}
\frac{\partial u}{\partial t}+\frac{\partial u}{\partial a}=D(a)\frac{%
\partial ^2{u}}{\partial x^2}-d(a)u, & t\geq t_0\geq A_l,\;a\geq 0,\; 0< x< \pi,\\
w(t,x)=\int_r^{A_l}u(t,a,x)da,& t\geq t_0\geq A_l,\; 0<x<\pi,\\
u(t,0,x)=f(w(t,x)), & t\geq t_0\geq A_l,\; 0<x<\pi,\\
u_x(t,a,0)=u_x(t,a,\pi)=0,\;& t\geq t_0\geq A_l,\; a \geq 0,
\end{array}
\right.
\end{equation*}
when the diffusion rate $D(a)$ and the death rate $d(a)$  are age dependent, and when the birth function $f(w)$ is nonmonotone. We also presented some illustrative examples. 
\end{abstract}


\vspace{0.2in} \noindent \textbf{Keywords and Phrases}: Steady-state
solution; Global Attractivity;  Age-structured model; Birth function; Age-dependent death rate; Age-dependent diffusion rate.

\noindent \textbf{AMS (1991) Subject Classification}: 34K10, 35K55, 35B35,
92D25


\section{Introduction}

Spatial movement and temporal maturation are two important characters in
most of biological systems; modeling the interaction between them has
attracted considerable attention recently \cite{Kuang,Gourley5,Gourley4,Gourley3,Huddleston,Liang,Mei,Ou4,Ou5,Ou3,Ou2,Smith,So1,Thieme1,Xu,Zhao2}%
. One of the most important methods applied is the Smith-Thieme
age-structure technique \cite{Smith}. In this approach, species population
is divided into two groups: mature and immature. At different ages, the
standard model with age-structured and diffusion is incorporated ( see \cite
{Metz}):
\begin{equation}
\frac{\partial u}{\partial t}+\frac{\partial u}{\partial a}=D(a)\frac{%
\partial ^2{u}}{\partial x^2}-d(a)u.  \label{equI1}
\end{equation}
Here $u:=u(t,a,x)$ denotes to the density of the population of the species at time $%
t\geq 0$, age $a\geq 0$, and location $x\in [0,\pi]$. The age functions $D(a)$ and $d(a)$ are the diffusion and death rates, respectively. Let $r\geq 0$
be the maturation time for the species and $A_l>0$ be the life span of the
species. Then the density of the mature population at time $t\geq 0$ and location $x\in [0,\pi]$ is given by
\begin{equation}
w(t,x)=\int_r^{A_l}u(t,a,x)da.  \label{equI2}
\end{equation}
Since only the mature individuals can reproduce, one can assume
\begin{equation}
u(t,0,x)=f(w(t,x))\,,  \label{equI3}
\end{equation}
where $f(.)$ is a birth function.

In \cite{So1}, So, Wu, and Zou assume that the diffusion and death rates, $%
D(a)$ and $d(a)$, of the mature population are age independent. i.e.,
\[
D(a)=D_m\;\;\;\text{and\ }\;\;d(a)=d_m.
\]
By this assumption and by substituting (\ref{equI1}) into (\ref{equI2}), they 
derived the following reaction diffusion equation:
\begin{equation}
\frac{\partial w}{\partial t}=D_m\frac{\partial ^2w}{\partial x^2}%
-d_mw+u(t,r ,x),  \label{ee}
\end{equation}
where $u(t,r ,x)$ is called the maturation rate and it can be obtained by the Fourier transforms from
(\ref{equI1}) and the boundary condition (\ref{equI3}), with a formula given
by
\[
u(t,r ,x)=\epsilon \int_{-\infty }^\infty f(w(t-r,y))K_\alpha (x-y)dy,
\]
where
\[
\epsilon :=\exp \left[ -\int_0^rd_I(a)da\right] ,\;\;\alpha :=\int_0^rD_I(a)da,
\]
and
\[
K_\alpha (x)=\frac{\exp ^{\left( -x^2/4\alpha \right) }}{\sqrt{4\pi \alpha }}%
.
\]
The functions $D_I(a)$ and $d_I(a)$ given above are the diffusion and death rates of the immature population, respectively.

As such, a non-local time-delayed reaction diffusion equation for the mature population can be obtained:
\[
\frac{\partial w}{\partial t}=D_m\frac{\partial ^2w}{\partial x^2}%
-d_mw+\epsilon \int_{-\infty }^\infty b(w(t-r,y))f_\alpha (x-y)dy.
\]
During the past decade, there have been some further studies on
this model. In \cite{Mei}, {Mei and So} investigated the stability of
traveling wave solution in the case of Nicholson's blowflies birth function.
Liang and Wu in \cite{Liang4} investigated the existence of traveling wave
solutions for different birth functions. In \cite{Liang}, Liang, So, Zhang,
and Zou considered the above model on a bounded domain where they assumed the diffusion
and death rates of the mature population to be constants. In fact, they  investigated the long time behavior of the solution by using a numerical simulation.

Thieme and Zhao in \cite{Thieme1} considered the following general
stage-structured model:
\begin{equation}\label{XuZhaoint}
\left\{
\begin{array}{ll}
\partial _tu+\partial _au=d_I(a)\Delta _xu-\mu _I(a)u, & 0<a<r,\;\;x\in %
\mathbb{R}^n \\
u(t,0,x)=f(u_m(t,x)), & \;\;t\geq -r,x\in \mathbb{R}^n, \\
\partial _tu_m=D_m\Delta _xu_m-d_m g(u_m)+u(t,r,x), & \;t>0,x\in \mathbb{R}^n,
\end{array}
\right.
\end{equation}
where $u_m$ and $u$ are the population density of the mature and the immature populations, $f(u_m)$ and $g(u_m)$ are the birth and death functions, $D_I(a)$ and $%
\mu_I(a)$ are the diffusion and death rates of the immature population, and $%
D_m$ and $d_m$ are age independent diffusion and death rates of the mature
population. In fact, they investigated the existence of traveling wave solutions of
this model, when the spatial domain is  $\R^n$.
When the spatial domain is a bounded region $\Omega\subset \R^n$, the  model given in Eq.  \eqref{XuZhaoint} was investigated by Xu and
Zhao in \cite{Xu}, and by Jin and Zhao in \cite{Zhao2}. In fact, the authors 
investigated the existence and the global attractivity of the steady-state solutions when the function $f(u_m)$ is monotone. In these articles the authors  assumed the diffusion and death rates of the mature population to be age-independent. i.e.,  $D(a)=D_m$ and $d(a)=d_m$, respectively.  Under these assumptions they transformed \eqref{XuZhaoint} into the following non-local time-delayed reaction diffusion equation:
\begin{equation}\label{reactiondiffusionint}
\left\{
\begin{array}{ll}
\partial _tu_m+\partial _au_m=D_m\Delta _xu_m-d_mg(u_m)+\\ \displaystyle \int_\Omega\Gamma(\eta(\tau),x,y)\mathcal F(\tau)f(u_m(t-\tau,y))dy, &  t>0,\;\;x\in %
\Omega, \\
Bu_m(t,x)=0, & \;\;t\geq 0,x\in \partial\Omega,  \\
u_m(t,x)=\phi(t,x), & \;t\in[-\tau,0],x\in \Omega,
\end{array}
\right.
\end{equation}
where $\Omega$ is a bounded region in $\R^n$, $\Gamma(\eta(\tau),x,y)$ is the Green's function associated with the Laplacian operator $\Delta_x$, $Bu_m=\frac{\partial u_m}{\partial n}+\alpha u$, $\eta(a)=\int_0^aD_I(s)ds$, $\mathcal F(a)=e^{-\int_0^a \mu_i(s)ds}$, and $\phi(t,x)$ is positive initial function. As a special case of the this model, Zhao in  \cite{Zhaoquartarly}  considered the following 
time-delayed reaction diffusion equation:
\begin{equation}\label{reactiondiffusionint2}
\left\{
\begin{array}{ll}
\displaystyle\frac{\partial u(t,x)}{\partial t}=D\Delta u(t,x)-\alpha u(t,x)+\\ \displaystyle \int_\Omega K(x,y) f(u(t-\tau,y))dy, &  t>0,\;\;x\in %
\Omega, \\\\
  
\displaystyle \frac{\partial u}{\partial n}=0, & \;\;t\geq 0,x\in \partial\Omega,   \\
u(t,x)=\phi(t,x), & \;t\in[-\tau,0],x\in \Omega,
\end{array}
\right.
\end{equation}
where $\Omega$ is bounded region in $\R^n$, $D,\alpha >0$, $\tau\geq 0$, $\Delta$ is the Laplacian operator, and $\frac{\partial u}{\partial n}$ is the normal derivative of $u$ in the direction of the outer normal $n$ to the $\partial \Omega.$ In fact, the author proved the global attractivity of the positive  steady state when $f(u_m)$ is a nonmonotone function in $u_m$. Conclusively, all these studies assumed the diffusion and death rates of the mature populations to be age-independent. 

Biologically, it is more realistic to include the age effects in the mathematical models during the whole life of the species. For
example,  women in the age between 15-40 years have higher birth
rate and lower death rate. This causes a variation in the diffusion and death
rates among the different ages of the mature individuals. Therefore, the authors of \cite{AlJararha} investigated the age-structured model (\ref{equI1})--(\ref{equI3}) when the diffusion and death rates are age-dependent. For the spatial domain they considered two cases. The first case is when the spatial domain is whole r4eal line $\R$. In this case they investigated the existence of monotone traveling waves solutions. The second case is when the spatial domain is closed and bounded interval in $\R$. For this case  they considered the model (\ref{equI1})--(\ref{equI3}) with different types of boundary conditions. Particularly, for the Neumann boundary conditions, $u_x(t,a,0)=u_x(t,a,\pi)=0$, the authors derived the following integral equation:
\begin{equation}\label{equ18}
w(t,x)=\displaystyle\int_r^{A_l}\int_0^\pi f(w(t-a,y))K(a,x,y)dyda,\;\;t\geq A_l.  
\end{equation}
The kernel function  $K(a,x,y)$  is given by 
\begin{equation*}
K(a,x,y) =\frac{\beta (a)}\pi \left( 1+\displaystyle %
2\sum_{n=1}^\infty e^{-n^2\alpha (a)}\cos nx\,\cos ny\right),
\end{equation*}
where $\alpha(a):=\int_0^aD(\xi)d\xi$ and $\beta(a):=e^{-\int_0^ad(\xi)d\xi}.$
For this integral equation and under certain conditions on the birth function $f(w)$, the authors of \cite{AlJararha} investigated the existence  of the positive  steady state solution $w^*$. Obviously, $w^*$ solves the algebraic equation $w^*=k^*f(w^*)$, where $k^*:=\int_r^{A_l}\beta(a)da$. Moreover, they investigated the global attractivity of the  positive steady state solution $w^*$ when the birth function $f(w)$ is  monotone. 
In this paper, we  investigate the global attractivity of the  positive and  steady state solution $w^*$ when the birth function $f(w)$ is nonmonotone. 
The paper is organized as follows. In section 2, we  present some preliminary results. In section 3, we prove our main result. In section 4, we present some illustrative examples. 
\section{\textbf{Preliminaries}}
In this section, we present some preliminary  results. Let $\mathbb X=C([0,\pi])$ be the space of continuous functions on the closed interval $[0,\pi]$ with the supremum norm defined by 
  $\|\phi\|_\infty=\sup_{x\in[0,\pi]}\phi(x)$.  Let $\mathbb X_+=\left\{\phi(x)\in \mathbb X\;| \;\;\phi(x)\geq 0\right\}$ be its positive cone. Then $\mathbb X_+$ has a non-empty interior. Hence, we can define a strongly positive relation on $\mathbb X_+$. In fact, for $\phi,\;\psi \in \mathbb X$, we have $\phi(x)\preceq \psi(x)$ if and only if  $\phi(x)\leq \psi(x),\;\forall x \in [0,\pi]$. Hence, ($\mathbb X, \mathbb X_+)$ is strongly ordered Banach space. We also consider the space $\mathbb Y=C([t_0-A_l, t_0], \mathbb X)$, where $t_0\geq A_l$ is fixed, with its ordered positive cone $\mathbb Y_+=C([t_0-A_l,t_0],\mathbb X_+)$. For convenience, we identify each $\phi \in \mathbb Y_+$ as a function from $[t_0-A_l, t_0]\times [0,\pi]$ to $\mathbb{R}$ as follows: $\phi(s,x)=\phi(s)(x)$. For any function $y(.):\left[t_0-A_l,b\right)\rightarrow \mathbb X$, where $b>t_0$, define $y_t \in \mathbb Y$ by $y_t(s)=y(t+s)$, $\forall s\in [t_0-A_l,t_0)$. For given $\beta >0$, we define the positive cone  $\Sigma_\beta:=\left\{\phi(x)\in \mathbb X_+\left| \;\;\phi(x)\leq\beta\right\}\right.$ and the function space  $\mathbb Z_\beta=C([t_0-A_l,t_0],\Sigma_\beta)$.
By using the method of steps, see, e.g., \cite{WuBook}), a unique solution $w(t,x,\phi )$ of the integral equation
\begin{equation}\label{semiflow_equ}\left\{
\begin{array}{lc}
w(t,x)=\displaystyle \int_r^{A_l}\int_0^\pi f(w(t-a,y))K%
(a,x,y)dyda, & \; t\geq A_l,\; x,y \in[0,\pi], \\
w_x(t,0)=w_x(t,\pi )=0, & t\geq t_0, \\
w(s,x)=\phi (s,x)\geq 0, & t_0-A_l\leq s\leq t_0,\;x\in [0,\pi], 
\end{array}
\right.
\end{equation}
where 
\begin{equation}\label{kernel}
K(a,x,y) =\frac{\beta (a)}\pi \left( 1+\displaystyle %
2\sum_{n=1}^\infty e^{-n^2\alpha (a)}\cos nx\,\cos ny\right),
\end{equation}
globally exists for any $\phi \in \mathbb Y_{+}$, provided that $f(w)$ is a Lipschitz function.
Therefore, we can define the semiflow $\Phi(t):\mathbb Y_{+}\rightarrow \mathbb Y_{+}$, by $(\Phi(t)\phi)(s,x)=w(t+s,x,\phi),\, \forall s\in [t_0-A_l,t_0],\, x\in[0,\pi]$. Moreover, the semiflow $\Phi
(t):\mathbb Y_{+}\rightarrow \mathbb Y_{+}$ is compact for $\forall t> t_0$. The concept of the semiflow can be found  in  \cite[p.8]{ZhaoBook} ( one can also see   \cite[p. 2]{Smith1}). We note that the kernel function in Eq. \eqref{kernel} is continuous and positive, and it uniformly bounded on $ [r,A_l]\times[0,\pi]\times[0,\pi]$, see, \cite{AlJararha}). Moreover,  $\int_r^{A_l}\int_0^\pi K
(a,x,y)dyda=k^*,$ where $k^*:=\int_r^{A_l}\beta(a)da$. 
To prove our main result, we need the following assumptions on the birth function $f(w)$: 
\begin{itemize}
\item [(F)] Assume that
\begin{itemize}
\item[(F1)] $f:\R^+\rightarrow \R$ is  Lipschitz continuous function $\forall w\geq 0$, $f(0)=0$, $f$ is differentiable at $0$ and $f^{\prime}(0)=p>0,\; \text{and}\; f(w)\leq p w, \forall w\geq0 .$
\item[(F2)] there exists a positive constant $M$, such that $\forall w>M$ we have $k^*\bar{f}(w)\leq w$, where $\bar f(w):=\max_{v\in [0,w]}f(v).$
\end{itemize}
\end{itemize}
Now, consider the following linearization of \eqref{semiflow_equ}:
\begin{equation}\label{linearization_equ}\left\{
\begin{array}{lc}
w(t,x)=\displaystyle p\int_r^{A_l}\int_0^\pi w(t-a,y)K%
(a,x,y)dyda, & \; t\geq A_l,\; x,y \in[0,\pi], \\
w_x(t,0)=w_x(t,\pi )=0, & t\geq t_0, \\
w(s,x)=\phi (s,x)\geq 0, & t_0-A_l\leq s\leq t_0,\;x\in [0,\pi], 
\end{array}
\right.
\end{equation}
Let $w(t,x)=e^{\lambda t}w(x)$ in the above equation, we get the following eigenvalue problem:
\begin{equation}\label{eigen_equ}\left\{
\begin{array}{lc}
w(x)=\displaystyle p\int_r^{A_l}\int_0^\pi e^{-a\lambda }w(y)K%
(a,x,y)dyda, & \; x \in[0,\pi], \\
w_x(0)=w_x(\pi )=0. &  
\end{array}
\right.
\end{equation}
Let $w(x)=1$. Then the characteristic equation is given by $p\Gamma_0(\lambda)=1$, where $\Gamma_0(\lambda):=\int_r^{A_l}\int_0^\pi e^{-a\lambda }K(a,x,y)dyda$. By solving this characteristic equation, we can determine uniquely the principle eigenvalue. For the kernel function given in equation \eqref{kernel}, we have $\Gamma _0(\lambda )=\int_r^{A_l}\exp\left\{-\left( \lambda \,a+\gamma (a)\right)\right\}da$ and $\gamma (a):=\int_0^ad(\xi )d\xi$. Clearly, $\Gamma_0(\lambda)$ is decreasing function in $\lambda$, and  $$e^{( -\lambda \,r)} \int_r^{A_l}e^{-\gamma (a)}da\leq \int_r^{A_l}\exp\left\{-\left( \lambda \,a+\gamma (a)\right)da\right\}\leq e^{( -\lambda \,A_l)} \int_r^{A_l}e^{-\gamma (a)}da.$$ Hence, $\displaystyle\lim_{\lambda\rightarrow \infty}\Gamma_0(\lambda)=0$,  $\displaystyle\lim_{\lambda\rightarrow -\infty}\Gamma_0(\lambda)=\infty$, and $\Gamma_0(0)=k^*:=\int_r^{A_l}e^{-\gamma (a)}da$. Therefore,  $p\Gamma_0(\lambda)=1$ has a unique solution $\lambda_0$ which represents the principle eigenvalue of \eqref{eigen_equ} \cite[Theorem 5.1]{AlJararha}. Moreover, $\lambda_0>0$ if $pk^*>1$ and $\lambda_0<0$ if $pk^*<1$. 

By applying the same argument used in the proof of Lemma 6.1 in \cite{AlJararha}, we have the following theorem:
\begin{theorem}
\label{connattlemma}
Assume that (F1) and (F2) hold. Then for any $\phi \in \mathbb Y_{+}$, a unique solution $w(t,x,\phi )$ of \eqref{semiflow_equ} exists, and $\lim\sup_{t\rightarrow \infty }w(t,x,\phi)\leq M$ uniformly $\forall x\in [0,\pi]$. Furthermore, the semiflow $\Phi (t):\mathbb Y_{+}\rightarrow \mathbb Y_{+}$ admits
a connected global attractor on $\mathbb Y_{+}$ which attracts every bounded set in $\mathbb Y_+$.
\end{theorem} 
Also, by applying the same argument used in the proof of Lemma 6.2 and Theorem 6.3 in \cite{AlJararha}, we have the following theorem:
\begin{theorem}
\label{unformpersts-andzeroconv}
Let F1 and F2 hold, and let $w(t,x,\phi)$ be a solution of \eqref{semiflow_equ} for $\phi\in \mathbb Y$. Then the following statements are valid:
\begin{itemize}
\item [I)] If $pk^*<1$  and $\phi\in \mathbb Y$, then $\displaystyle\lim_{t\rightarrow\infty} w(t,x,\phi)=0.$ 
\item [II)] If $pk^*>1$ , then \eqref{semiflow_equ} admits at least one homogeneous steady state solution $w^*\in[0,M]$, and there exists a positive constant $\delta$ such that $\displaystyle \liminf_{t\rightarrow \infty }w(t,x,\phi)\geq \delta$ uniformly, for all $\phi \in \mathbb Y_{+}$ and $x\in [0,\pi]$.
\end{itemize}
\end{theorem} 
\begin{remark} Assume that $pk^*>1$, and (F1) and (F2) hold. Let $F(w)=k^*f(w)-w$. Since $f(w)$ satisfies (F1). Then $F(0)=0$ and $F^\prime(0)=pk^*-1>0$. Moreover, since $f(w)$ satisfies  (F2), then $F(M)\leq 0$. Therefore, there exists some $w^*\in (0,M]$ such that $F(w^*)=0.$ Hence, $w^*$ is a positive constant steady state of \eqref{semiflow_equ}.
\end{remark} 
\section{The Main Result} 
In this section, we prove the global attractivity of the positive  steady state solution $w^*$. To prove this result, we use the fluctuation method, e.g., see, \cite{Hsu,Thieme,Thieme1,Thieme2,Thieme3,Zhaoquartarly}. To apply this method, we need the following assumption on the birth function $f(w)$:
\begin{itemize}
\item [(F3)] $f^\prime(0)>1,\;\text{and}\;\frac{f(w)}{w}$ is strictly decreasing for $u \in (0,M]$, and $f(w)$ satisfies the property (P): that is for any $u,v\in (0,M]$ and $u\leq w^* \leq v,$ $ u\geq k^* f(v)$ and $v\leq k^* f(u)$, we have $v=u$.
\end{itemize}
\begin{lemma}(Lemma 2.2 \cite{Hsu}): A function $f(w)$ satisfies the property (P) if one of the following conditions hold:
\begin{itemize}
\item[(P0)] $f(w)$ is non-decreasing on $[0,M]$.
\item[(P1)] $wf(w)$ is strictly increasing on $(0,M]$.
\item[(P2)] $f(w)$ is non-increasing for $w\in [w^*,M]$, and $\frac{f(k^*f(w))}{w}$ is strictly decreasing for all $w\in (0,w^*]$.
\end{itemize}
\end{lemma}
\begin{lemma}\label{mainlem}
Let $\phi\in \mathbb Y_+$ with $\phi(t_0,.)\not\equiv 0.$ Moreover, Let $\omega(\phi)$ be the omega limit set of the positive orbits through $\phi$ for the solution semiflow $\Phi(t).$  Then $\mathbb Z_M$ is positively invariant. i.e., $\Phi(t)\mathbb Z_M\subset \mathbb Z_M$. moreover, $\omega(\phi)\subset  \mathbb Z_M.$ 
\end{lemma}
\noindent\emph{Proof.} Let $\phi\in \mathbb Y_+$ with $\phi(t_0,.)\not\equiv 0.$ Moreover, Let $\omega(\phi)$ be the omega limit set of the positive orbits through $\phi$ for the solution semiflow $\Phi(t).$ From Theorem \ref{connattlemma}, we have $\lim\sup_{t\rightarrow \infty }w(t,x,\phi)\leq M$, $\forall x \in [0,\pi]$. Therefore, $\Phi(t)\mathbb Z_M\subset \mathbb Z_M$, and hence, $\omega(\phi)\subset \mathbb Z_M.$ $\square$

\begin{theorem}\label{maintheorem}
Assume that $pk^*>1$. Moreover, assume that (F1) and (F2) hold. Then for any $\phi \in \mathbb Y_+$ with  $\phi(t_0,.)\not\equiv 0,$ we have $\displaystyle \lim_{t\rightarrow \infty} w(t,x,\phi)=w^*$ uniformly $\forall x \in [0,\pi].$
\end{theorem}
\noindent{\emph{Proof.}} To prove the global attractivity of $w^*$, by Lemma \ref{mainlem},  it is sufficient to prove the globall attractivity of $w^*$ on  $\mathbb Z_M$. Therefore, let $\phi \in \mathbb Z_M$ be such that $\phi(t_0,.)\not\equiv 0.$ Then the solution of \eqref{semiflow_equ} through $\phi$;  $w(t,x,\phi):=w(t,x)$ satisfies 
\[
w(t,x)=\int_r^{A_l}\int_0^\pi K(a,x,y)f(w(t-a,y))dyda.
\]
Let $w^\infty(x)=\displaystyle \limsup_{t\rightarrow \,\infty}w(t,x)$ and $w_\infty(x)=\displaystyle \liminf_{t\rightarrow \infty}\,w(t,x)$ for any $x\in [0,\pi].$ Then $w^\infty(x)\geq w_\infty(x)$. Since $pk^*>1$. Then by Theorem   \ref{unformpersts-andzeroconv}, we have 
\[
0<\delta\leq w_\infty(x)\leq w^\infty(x)\leq M.
\]
Moreover, if we let $w^\infty=\displaystyle\sup_{x\in[0,\pi]}w^\infty(x)$ and $w_\infty=\displaystyle\inf_{x\in[0,\pi]}w_\infty(x)$, then $0<\delta\leq w_\infty\leq w^\infty\leq M.$ Now, define the  function
\[
F(u,v)=\left\{ \begin{array}{ll}
\min \{f(w):\; u\leq w\leq v\},\;& \text{if} \;u\leq v,\\
\max \{f(w):\;v\leq w\leq u\},\;& \text{if} \;v\leq u. 
\end{array}\right.
\]
Then $F(u,v):[0,M]\times[0,M]\rightarrow \R$ is continuous function, nondecreasing in $u\in [0,M]$, nonincreasing in $v\in [0,M]$, and $f(w)=F(w,w)$, see, e.g., \cite[section  3.6]{Thieme2}.
Since the kernel function $K(a,x,y)$ is uniformly bounded $\forall (a,x,y)\in[r,A_l]\times[0,\pi]\times[0,\pi]$. Then, by Fatou's lemma, we have
\begin{eqnarray*}
w^\infty(x)=\displaystyle \limsup_{t\rightarrow \infty}w(t,x)&=&\displaystyle \limsup_{t\rightarrow \infty}\int_r^{A_l}\int_0^\pi K(a,x,y)f(w(t-a,y))dyda\\
&\leq&\int_r^{A_l}\int_0^\pi K(a,x,y)\limsup_{t\rightarrow \infty}f(w(t-a,y))dyda\\
&=&\int_r^{A_l}\int_0^\pi K(a,x,y)\limsup_{t\rightarrow \infty}F(w(t-a,y),w(t-a,y))dyda\\
&\leq& \int_r^{A_l}\int_0^\pi K(a,x,y)F(w^\infty ,w_\infty)dyda\\
&=&k^*F(w^\infty ,w_\infty).
\end{eqnarray*}
Hence, 
\begin{equation}\label{inq1}
w^\infty(x)\leq k^*F(w^\infty,w_\infty). 
\end{equation}
Similarly, we have 
\begin{equation}\label{inq2}
w_\infty(x)\geq k^*f(w_\infty,w^\infty).
\end{equation} 
Clearly, by the definition of $F(u,v)$ there exists $u,v\in[w_\infty,w^\infty]\subset [0,M]$ such that $f(u)=F(w^\infty,w_\infty)$ and $f(v)=F(w_\infty,w^\infty)$. Hence, 
\begin{equation}\label{inq3}
f(u)=F(w^\infty,w_\infty)\geq \frac{w^\infty}{k^*}\geq\frac{u}{k^*}\,\left(\frac{v}{k^*}\right),
\end{equation}
and 
\begin{equation}\label{inq4}
f(v)=F(w_\infty,w^\infty)\leq \frac{w_\infty}{k^*}\leq\frac{v}{k^*}\,\left(\frac{u}{k^*}\right).
\end{equation}
Therefore,
\[
\frac{k^*f(v)}{v}\leq 1=\frac{k^*f(w^*)}{w^*}\leq \frac{k^*f(u)}{u}.
\]
Since $\frac{f(w)}{w}$ is assumed to be a strictly decreasing function on $(0,M]$. Then $u\leq w^* \leq v.$ Also, by \eqref{inq3} and \eqref{inq4}, we have 
\begin{equation}
f(u)\geq \frac{w^\infty}{k^*}\geq\frac{v}{k^*},
\end{equation}
and 
\begin{equation}\label{inq6}
f(v)\leq \frac{w_\infty}{k^*}\leq\frac{u}{k^*}.
\end{equation}
i.e., 
\[
k^*f(u)\geq w^\infty\geq v\;\;\; \text{and}\;\;\; k^*f(v)\leq w_\infty\leq u.
\]
Since $f(w)$ satisfies th property (P), then we get  $w^*=u=v.$  Moreover, we have 
\[
k^*f(u)\geq w^\infty\geq u \;\;\; \text{and} \;\;\;f(v)\leq w_\infty\leq v.
\]
Hence, $w^*=w_\infty=w^\infty.$ Recall that 
\[
w^\infty \geq w^\infty(x)\geq w_\infty(x)\geq w_\infty,\;\;\; \forall x \in [0,\pi].
\]
Thus, $w^\infty(x)= w_\infty(x)=w^*, \; \forall x \in [0,\pi].$ Hence,
\begin{equation}\label{limt}
\displaystyle\lim_{t\rightarrow \infty} w(t,x)=w^*,\; \forall x \in [0,\pi].
\end{equation}
Finally, we show that $\displaystyle\lim_{t\rightarrow \infty} w(t,x)=w^*$ uniformly $\; \forall x \in [0,\pi].$ In fact, it is enough to show that $\omega(\phi)=\{w^*\},\forall \phi \in \mathbb Y_+.$ Let $\psi \in \omega(\phi).$ Then there exists a time sequence $t_n\rightarrow \infty$ such that $\Phi(t_n)\phi\rightarrow \psi$ in $\mathbb Y$ as $n\rightarrow \infty.$ This implies that
\[
\displaystyle \lim_{n\rightarrow \infty}w(t_n+s,x,\phi)=\psi(s,x)
\]
uniformly for $(s,x)\in [t_0-A_l,t_0]\times[0,\pi].$ Hence, from \eqref{limt}, we have $\psi(s,x)=w^*, \forall (s,x)\in [t_0-A_l,t_0]\times[0,\pi].$ Thus, we get $\omega(\phi)=w^*$, which implies that $w(t,.,\phi)$ converges to $w^*$ in $\mathbb X$ as $t\rightarrow \infty.$
\section{Examples}
In this section, we present three examples of the birth function $f(w)$ to demonstrate the application of our main result. 

First, we begin with the Ricker type function $f(w)=pwe^{-aw^{q}}$,  $a$,  $p>0,$ and $q>0$. Then, we have the following theorem: 
\begin{theorem}
Let $f(w)=pwe^{-aw^q},$ $a>0$,  $p>0$, $q>0$.  Moreover,  assume that  $1< k^* p\leq e^{\frac{2}{q}}$. Then the unique positive  steady state solution $w^*=\left[\frac{1}{a}\ln(pk^*)\right]^{\frac{1}{q}}$ attracts all positive solutions of \eqref{semiflow_equ}.
\end{theorem} 
\noindent \emph{Proof.} First, we remark that $f(w)$ satisfies the conditions (F1)-(F3), $f(w)/w$ is a strictly decreasing function on $[0,\infty)$. Moreover, $f^{\prime}(0)=p>0$, and $f(w)$ takes its maximum at $\overline{w}=(\frac{1}{aq})^{\frac{1}{q}}$ and $f(\overline{w})=p(\frac{1}{aqe})^{\frac{1}{q}}$. Assume that $1<k^*p\leq e^{\frac{1}{q}}$, then $f(w)$ is a monotone increasing function on $[0,w^*]$. Therefore, (P0) holds with $M=w^*$. Now,  assume that $k^*p>e^{\frac{1}{q}}$. In this case, we consider $M=f(\overline{w})$. Hence, $f(w)$ is decreasing function on $[w^*,M]$. Moreover, the function
\[
h(w):=\frac{f(k^*f(w))}{w}=p^2k^*\exp\left\{-a\left(w^q+(pk^* w)^q e^{-aqw^q} \right) \right\}.
\] 
 is a strictly decreasing function on $[0,u^*]$ if $e^{\frac{1}{q}}<pk^* \leq e^\frac{2}{q}$. Therefore, the property (P2) holds. Hence, the conditions of Theorem  \ref{maintheorem} hold. Thus, $w^*$ attracts every positive solution of \eqref{semiflow_equ}.$\square$
 \\
 
Next, we consider the Beverton-Holt function $f(w)=\frac{pw}{1+aw^q}$, $a>0$, $p>0$, and $q>0.$ Then we have the following theorem:
\begin{theorem}
Let $f(w)=\frac{pw}{1+aw^q}$, $a>0$, $p>0$,  $q>0$.  Moreover,  assume that $q\in \left(0,\max \left(2,\frac{pk^*}{pk^*-1}\right)\right]$, or $q>\max \left(2,\frac{pk^*}{pk^*-1}\right)$ and $k^*f(\overline w)\leq \left(\frac{2}{a(q-2)} \right)^{\frac{1}{q}}$; where $pk^*>1$ and $\overline w$ is the value where $f(w)$ takes its maximum. Then the unique positive  steady state solution $w^{*}=\left(\frac{pk^*-1}{a}\right) ^{1/q}$ attracts all positive solutions of \eqref{semiflow_equ}. 
\end{theorem}
\noindent\emph{Proof.} First, we remark that $f(w)$ satisfies the conditions (F1)-(F3), and $f(w)/w$ is a strictly decreasing function on $[0,\infty)$. Moreover, $f^{\prime}(0)=p>0$,  and $f(w)$ takes its maximum at $\overline{w}=(\frac{1}{a(q-1)})^{\frac{1}{q}}$ and $f(\overline{w})=\frac{p(q-1)}{q} \overline w$. Assume that $q \in (0,1]$, then $f(w)$ is monotone increasing on $[0,\infty)$, and hence, (P0) holds with $M=w^*$. Now, if we assume $1<q\leq 2$, then $wf(w)$ is increasing function on $[0,\infty)$. Hence (P1) holds with $M=w^*$. Moreover, if $1<pk^* \leq \frac{q}{q-1}$ (i.e., $q \in (1,\frac{pk^*}{pk^*-1})$), then  $w^*\leq \overline w$. Hence, if we let $M=w^*$, then (P0) holds. conclusively, if $q\in \left(0,\max \left(2,\frac{pk^*}{pk^*-1}\right)\right]$, then either (P0) or (P1) holds. If $q>\max \left(2,\frac{pk^*}{pk^*-1}\right)$, then $h(w):=wf(w)=\frac{p w^2}{1+a w^q}$ is a monotone increasing function on $\left[0, \left(\frac{2}{a(q-2)} \right)^{\frac{1}{q}} \right]$. Hence, if we consider $M=k^*f(\overline w)$, then (P1) holds provided that $k^*f(\overline w)\leq \left(\frac{2}{a(q-2)} \right)^{\frac{1}{q}}$. Hence, the conditions of Theorem  \ref{maintheorem} hold. Therefore, $w^*$ attracts every positive solution of \eqref{semiflow_equ}.$\square$
\\

Finally, we consider the logistic type birth function $f(w)=pw(1-\frac{w}{K})$, $p>0$, and $K>0$. Then we have the following theorem:
\begin{theorem}
Let $f(w)=pw(1-\frac{w}{K})$, $p>0$, and $K>0$ in \eqref{semiflow_equ}. Moreover, assume that  $1<pk^*\leq 3$.  Then the unique positive steady state solution $w^*=K\left(1-\frac{1}{pk^*}\right)$ attracts every positive solution of \eqref{semiflow_equ}.
\end{theorem} 
\noindent\emph{Proof.}
First, we remark that $f(w)$ satisfies the conditions (F1)-(F3), and $f(w)/w$ is a strictly decreasing function on $(0,K]$. Moreover, $f^{\prime}(0)=p>0$, and $f(w)$ takes its maximum at $\overline{w}=\frac{K}{2}$ with $f(\overline{w})=\frac{pK}{4}$. Assume that $1<pk^*\leq 2$, then $f(w)$ is a monotone increasing function on $[0,\frac{2}{K}]$, and hence, (P0) holds with $M=u^*$. Assume that $2<pk^*<4$, and let $M=\frac{pk^*}{4}K$. Then the function
\[
h(w):=\frac{f(k^*f(w))}{w}=\frac{p^2k^*}{K^3}\left(K^2\left(K-w\right)-pk^*w\left(K-w\right)^2\right). 
\] 
is a strictly decreasing function on $[0,w^*]$ provided that $2<pk^*\leq 3$. Hence, the property (P2) holds. Therefore, the assumptions of Theorem \ref{maintheorem} hold. Thus $w^*$ attracts every positive solution of \eqref{semiflow_equ}.$\square$

\end{document}